\documentstyle[12pt,twoside]{article}

\font\teneufm=eufm10 scaled \magstep1
\font\seveneufm=eufm7 scaled \magstep1
\font\fiveeufm=eufm5  scaled \magstep1
\newfam\eufmfam
\textfont\eufmfam=\teneufm
\scriptfont\eufmfam=\seveneufm
\scriptscriptfont\eufmfam=\fiveeufm

\newfam\msbfam
\font\tenmsb=msbm10 scaled \magstep1  \textfont\msbfam=\tenmsb
\font\sevenmsb=msbm7 scaled \magstep1 \scriptfont\msbfam=\sevenmsb
\font\fivemsb=msbm5 scaled \magstep1  \scriptscriptfont\msbfam=\fivemsb
\def\Bbb{\fam\msbfam \tenmsb}

\def\RR{{\Bbb R}}
\def\CC{{\Bbb C}}

\def\ra{\rightarrow}
\def\HollowBoxx #1#2#3{{\dimen0=#1 \advance\dimen0 by -#2
       \dimen1=#1 \advance\dimen1 by #3
        \vrule height 0pt depth #3 width #2
       \hskip -#3
       \vrule height #1 depth #3 width #3}}
 \def\LeftContraction{\mathord{\kern1.45pt \HollowBoxx{6pt}{3.5pt}{.4pt}}\,}
 \def\HollowBox #1#2#3{{\dimen0=#1 \advance\dimen0 by -#3
       \dimen1=#1 \advance\dimen1 by #3
        \vrule height #1 depth #3 width #3
        \vrule height 0pt depth #3 width #2
        \hskip -#3}}
 \def\RightContraction{\mathord{\, \HollowBox{6pt}{3.1pt}{.4pt}} \kern1.6pt}
\def\qed{{\hfill $\Box$}}
\newtheorem{theorem}{THEOREM}[section]
\newtheorem{corollary}[theorem]{Corollary}
\newtheorem{lemma}[theorem]{Lemma}

\begin{document}
\begin{center}

{\Large \bf On the Dimension of the\\ \medskip
Stability Group for a\\ \medskip Levi Non-Degenerate Hypersurface}
\vspace{1cm}

\normalsize V. V. Ezhov and A. V. Isaev
\end{center}

\begin{quotation} \small \sl We classify locally defined Levi non-degenerate non-spherical real-analytic hypersurfaces in complex space for which the dimension of the group of local CR-automorphisms has the second largest positive value.    
\end{quotation}

\pagestyle{myheadings}
\markboth{V. Ezhov and A. Isaev}{Dimension of the Stability Group}

\setcounter{section}{0}
\section{Introduction}
\setcounter{equation}{0}

Let $M$ be a real-analytic hypersurface in $\CC^{n+1}$ passing through the origin. Assume that the Levi form of $M$ at 0 is non-degenerate and has signature $(n-m,m)$ with $n\ge 2m$. Then in some local holomorphic coordinates $z=(z_1,\dots,z_n),\, w=u+iv$ in a neighborhood of the origin, $M$ can be written in the Chern-Moser normal form (see \cite{CM}), that is, given by an equation
$$
v=\langle z,z\rangle+\displaystyle\sum_{k,\overline{l}\ge 2}F_{k\overline{l}}(z,\overline{z},u),
$$
where $\langle z,z\rangle=\displaystyle\sum_{\alpha,\beta=1}^nh_{\alpha\beta}z_{\alpha}\overline{z_{\beta}}$ is a non-degenerate Hermitian form with signature
$(n-m,m)$, and $F_{k\overline{l}}(z,\overline{z},u)$ are polynomials of
degree $k$ in $z$ and $\overline{l}$ in $\overline{z}$  whose coefficients are
analytic functions of $u$ such that the following conditions hold
\begin{equation}
\begin{array}{lll}
\hbox{tr}\,F_{2\overline{2}}&\equiv&0,\\
\hbox{tr}^2\,F_{2\overline{3}}&\equiv&0,\\
\hbox{tr}^3\,F_{3\overline{3}}&\equiv&0.
\end{array}\label{formnfcond}
\end{equation}
Here the operator $\hbox{tr}$ is defined as 
$$
\displaystyle \hbox{tr}:=\sum_{\alpha,\beta=1}^n\hat h_{\alpha\beta}\frac{\partial^2 }{\partial z_{\alpha}\partial \overline{z}_{\beta}},
$$
where $(\hat h_{\alpha\beta})$ is the matrix inverse to $H:=(h_{\alpha\beta})$. 
Everywhere below we assume that $M$ is given in the normal form.

Let $\hbox{Aut}_0(M)$ denote the group of all local CR-automorphisms of $M$ defined near 0 and preserving 0. To avoid confusion with the term \lq\lq isotropy group of $M$ at 0\rq\rq\, usually reserved for global CR-automorphisms of $M$ preserving the origin, this group is often called the {\it stability group}\, of $M$ at 0. Every element $\varphi$ of $\hbox{Aut}_0(M)$ extends to a biholomorphic mapping defined in a neighborhood of the origin in $\CC^{n+1}$ and therefore can be written as
$$
\begin{array}{lll}
z&\mapsto& f_{\varphi}(z,w),\\
w&\mapsto& g_{\varphi}(z,w),
\end{array}
$$
where $f_{\varphi}$ and $g_{\varphi}$ are holomorphic. We equip $\hbox{Aut}_0(M)$ with the topology of uniform convergence of the partial derivatives of all orders of the component functions on a neighborhood of 0. The group $\hbox{Aut}_0(M)$ with this topology is a topological group.

It follows from \cite{CM} that every element $\varphi=(f_{\varphi},g_{\varphi})$ of $\hbox{Aut}_0(M)$ is uniquely determined by a set of parameters $(U_{\varphi}, a_{\varphi}, \lambda_{\varphi},\sigma_{\varphi},r_{\varphi})$, where $\sigma_{\varphi}=\pm 1$,  $U_{\varphi}$ is an $n\times n$-matrix such that $\langle U_{\varphi}z,U_{\varphi}z\rangle=\sigma_{\varphi}\langle z,z\rangle$ for all $z\in\CC^n$, $a_{\varphi}\in\CC^n$, $\lambda_{\varphi}>0$, $r_{\varphi}\in\RR$ (note that $\sigma_{\varphi}$ can be equal to $-1$ only for $n=2m$). These parameters are determined by the following relations
$$
\begin{array}{ll}
\displaystyle\frac{\partial f_{\varphi}}{\partial z}(0)=\lambda_{\varphi}U_{\varphi}, &
\displaystyle\frac{\partial f_{\varphi}}{\partial w}(0)=\lambda_{\varphi}U_{\varphi}a_{\varphi},\\
\vspace{0.2cm} &\\
\displaystyle\frac{\partial g_{\varphi}}{\partial w}(0)=\sigma_{\varphi}\lambda_{\varphi}^2, &
\hbox{Re}\,\displaystyle\frac{\partial^2 g_{\varphi}}{\partial^2 w}(0)=2\sigma_{\varphi}\lambda_{\varphi}^2r_{\varphi}.
\end{array}
$$
For results on the dependence of local CR-mappings on their jets in more general settings see \cite{BER1}, \cite{BER2}, \cite{Eb}, \cite{Z}.   

We assume that $M$ is {\it non-spherical at the origin}, i.e., that $M$ in a neighborhood of the origin is not CR-equivalent to an open subset of the hyperquadric  given by the equation $v=\langle z,z\rangle$. In this case for every element $\varphi=(f_{\varphi},g_{\varphi})$ of  $\hbox{Aut}_0(M)$ the parameters $a_{\varphi}, \lambda_{\varphi},\sigma_{\varphi},r_{\varphi}$ are uniquely determined by the matrix $U_{\varphi}$, and the mapping 
$$
\Phi:\, \hbox{Aut}_0(M)\ra GL_n(\CC),\qquad \Phi:\,\varphi\mapsto U_{\varphi}
$$
is a topological group isomorphism between $\hbox{Aut}_0(M)$ and $G_0:=\Phi(\hbox{Aut}_0(M))$ with $G_0$ being a real algebraic subgroup of $GL_n(\CC)$; in addition the mapping
\begin{equation}
\Lambda:\, G_0(M)\ra \RR_{+},\qquad \Lambda:\, U_{\varphi}\mapsto\lambda_{\varphi}\label{formlambda}
\end{equation}
is a Lie group homomorphism with the property $\Lambda(U_{\varphi})=1$ if all eigenvalues of $U_{\varphi}$ are unimodular, where $\RR_{+}$ is the group of positive real numbers with respect to multiplication (see \cite{CM}, \cite{B}, \cite{L1}, \cite{BV}, \cite{VK}). Since $G_0(M)$ is a closed subgroup of $GL_n(\CC)$, we can pull back its Lie group structure to $\hbox{Aut}_0(M)$ by means of $\Phi$ (note that the pulled back topology is identical to that of $\hbox{Aut}_0(M)$). Let $d_0(M)$ denote the dimension of $\hbox{Aut}_0(M)$. We are interested in characterizing hypersurfaces for which $d_0(M)$ is large. 

If $n>2m$, $G_0(M)$ is a closed subgroup of the pseudounitary group $U(n-m,m)$ of all matrices $U$ such that
$$
U^tH\overline{U}=H,
$$
(recall that $H$ is the matrix of the Hermitian form $\langle z,z\rangle$). The group $U(n,0)$ is the unitary group $U(n)$. If $n=2m$, $G_0$ is a closed subgroup of the group $U'(m,m)$ of all matrices $U$ such that
$$
U^tH\overline{U}=\pm H,
$$
that has two connected components. In particular, we always have $d_0(M)\le n^2$. If $d_0(M)=n^2$ and $n>2m$, then $G_0(M)=U(n-m,m)$. If $d_0(M)=n^2$ and $n=2m$, then we have either $G_0(M)=U(m,m)$, or $G_0(M)=U'(m,m)$.

We will say that the group $\hbox{Aut}_0(M)$ is {\it linearizable}, if in some coordinates every $\varphi\in\hbox{Aut}_0(M)$ can be written in the form
\begin{equation}
\begin{array}{lll}
z&\mapsto& \lambda Uz,\\
w&\mapsto& \sigma\lambda^2 w.
\end{array}\label{formlinearform}
\end{equation}
Clearly, in the above formula $U=U_{\varphi}$, $\lambda=\lambda_{\varphi}$, $\sigma=\sigma_{\varphi}$. The group $\hbox{Aut}_0(M)$ is known to be linearizable, for example, for $m=0$ (see \cite{KL}) and for $m=1$ (see \cite{Ezh1}, \cite{Ezh2}). If all elements of $\hbox{Aut}_0(M)$ in some coordinates have the form (\ref{formlinearform}), we say that $\hbox{Aut}_0(M)$ is {\it linear}\, in these coordinates. It is shown in Lemma 3 of \cite{Ezh3} that if $\hbox{Aut}_0(M)$ is linear in some coordinates, it is linear in some normal coordinates as well.    

We will first discuss the case when $d_0(M)$ takes the largest possible value, that is, when $d_0(M)=n^2$. Observe that in this case $\hbox{Aut}_0(M)$ is linearizable for any $m$. Indeed, if $d_0(M)=n^2$, the group $G_0(M)$ contains $U(n-m,m)$. Hence $G_0(M)$ contains the subgroup $Q:=\{e^{it}\cdot E_n,\,t\in\RR\}$, where $E_n$ is the $n\times n$ identity matrix. Let $\hat Q=\Phi^{-1}(Q)\subset \hbox{Aut}_0(M)$. The subgroup $\hat Q$ is compact, and the argument in \cite{KL} (see also \cite{VK}) yields that in some normal coordinates every element of $\hat Q$ can be written in the form (\ref{formlinearform}). For every $\varphi\in \hat Q$ we clearly have $\sigma_{\varphi}=1$. Further, since $Q$ is compact, there are no non-trivial homomorphisms from $Q$ into $\RR_{+}$, and therefore $\lambda_{\varphi}=1$ for every $\varphi\in \hat Q$. Hence, in these coordinates the function
$$
F(z,\overline{z},u):=\displaystyle\sum_{k,\overline{l}\ge 2}F_{k\overline{l}}(z,\overline{z},u)
$$
is invariant under all linear transformations from $Q$ and thus $F_{k\overline{l}}\equiv 0$, if $k\ne \overline{l}$. We will now show that  
$\hbox{Aut}_0(M)$ is linearizable. Since linearizability arguments of this kind will occur several times throughout the paper, we give some details on the linearizability of $\hbox{Aut}_0(M)$ for general hypersurfaces.

Suppose that $M$ is given in the Chern-Moser normal form near the origin. The main step in showing that $\hbox{Aut}_0(M)$ is linearizable is to prove that in some normal coordinates for every $\varphi\in\hbox{Aut}_0(M)$, we have $a_{\varphi}=0$. Indeed, if $a_{\varphi}=0$, it follows from \cite{CM} that $\varphi$ in the given coordinates is a fractional linear transformation that becomes linear if $r_{\varphi}=0$. It is shown in the proof of Proposition 3 of \cite{L2} that $a_{\varphi}=0$ implies that $r_{\varphi}=0$, provided $\lambda_{\varphi}=1$. Furthermore, if for every $\varphi\in \hbox{Aut}_0(M)$ we have $a_{\varphi}=0$ and there exists $\varphi_0\in\hbox{Aut}_0(M)$ with $\lambda_{\varphi_0}\ne 1$, the group $\hbox{Aut}_0(M)$ becomes linear after applying a transformation of the form
\begin{equation}
\begin{array}{lll}
z&\mapsto&\displaystyle\frac{z}{1+qw},\\
\vspace{0cm}&&\\
w&\mapsto&\displaystyle\frac{w}{1+qw},
\end{array}\label{fractionallintransfm}
\end{equation}
for some $q\in\RR$. 

To prove that $a_{\varphi}=0$ for a fixed $\varphi=(f_{\varphi}, g_{\varphi})$ in the given coordinates, we introduce weights as follows. Let each of $z_1,\dots,z_n$, $\overline{z_1},\dots,\overline{z_n}$ be of weight 1 and $u$ be of weight 2. Then we can write a weight decomposition for the function $F$ as follows
$$
F(z,\overline{z},u)=\sum_{j=\gamma}^{\infty}F_j,
$$
where $F_j$ is the component of $F$ of weight $j$, and $F_{\gamma}\not\equiv 0$. Next, since $\varphi$ is a local automorphism of $M$, we have
\begin{equation}
\hbox{Im}\, g_{\varphi}=\langle f_{\varphi},f_{\varphi}\rangle + F\left(f_{\varphi},\overline{f_{\varphi}},\hbox{Re}\, g_{\varphi}\right),\label{eq1}
\end{equation}
where we set $v=\langle z,z\rangle +F(z,\overline{z},u)$. Extracting all terms of weight $\gamma+1$ from identity (\ref{eq1}), we obtain the following identity (see \cite{B}, \cite{L1}, \cite{L2})
\begin{equation}
\begin{array}{l}
\hbox{Re}\,\left(i\tilde g_{\gamma+1}+2\langle \lambda_{\varphi}^{-1}U_{\varphi}^{-1}\tilde f_{\gamma},z\rangle\right)|_{v=\langle z,z\rangle} + T(F_{\gamma},a_{\varphi})=\\
F_{\gamma+1}(z,\overline{z},u)-\displaystyle\frac{1}{\lambda_{\varphi}^2}F_{\gamma+1}(\lambda_{\varphi}U_{\varphi}z,\overline{\lambda_{\varphi}U_{\varphi}z},\lambda_{\varphi}^2u).
\end{array}\label{eq3}
\end{equation}
Here $\left(\sum_{j=1}^{\infty} \tilde f_j, \sum_{j=1}^{\infty} \tilde g_j\right)$ is the weight decomposition for the map $(\tilde f,\tilde g):=(f_{\varphi}-f_{\varphi}^Q,g_{\varphi}-g_{\varphi}^Q)$, where $\varphi^Q=(f_{\varphi}^Q,g_{\varphi}^Q)$ is the following local automorphism of the hyperquadric given by the equation $v=\langle z,z\rangle$
$$
\begin{array}{lll}
z&\mapsto & \displaystyle\frac{\lambda_{\varphi}U_{\varphi}(z+a_{\varphi}w)}{1-2i\langle z,a\rangle-(r_{\varphi}+i\langle a,a\rangle)w},\\
w&\mapsto & \displaystyle\frac{\sigma_{\varphi}\lambda_{\varphi}^2w}{1-2i\langle z,a\rangle-(r_{\varphi}+i\langle a,a\rangle)w},
\end{array}
$$
and
$$
\begin{array}{l}
\displaystyle T(F_{\gamma},a_{\varphi}):=2\hbox{Re}\,\Bigl(-2i\langle z,a_{\varphi}\rangle F_{\gamma}+(u+i\langle z,z\rangle)\sum_{j=1}^na_j\frac{\partial F_{\gamma}}{\partial z_j}+\\\displaystyle
2i\langle z,a\rangle \sum_{j=1}^n z_j\frac{\partial F_{\gamma}}{\partial z_j}+
i\langle z,a_{\varphi}\rangle(u+i\langle z,z\rangle)\frac{\partial F_{\gamma}}{\partial u}\Bigr),
\end{array}
$$
where $a_1,\dots,a_n$ denote the components of the vector $a_{\varphi}$.
 
If $F_{\gamma+1}=0$, the right-hand side of (\ref{eq3}) vanishes, and the proof of Proposition 1 in \cite{L1} shows that the resulting homogeneous identity can only hold if $a_{\varphi}=0$. Clearly, if $F_{k\overline{l}}\equiv 0$ for $k\ne \overline{l}$, then the weight decomposition for $F$ contains only terms of even weights, and, in particular, we have $F_{\gamma+1}=0$. Thus, we have shown that $\hbox{Aut}_0(M)$ is linearizable if $d_0(M)=n^2$.

Observe further that if $d_0(M)=n^2$, the mapping $\Lambda$ defined in (\ref{formlambda}) is constant, that is, $\lambda_{\varphi}=1$ for all $\varphi\in\hbox{Aut}_0(M)$. Indeed, consider the restriction of $\Lambda$ to $U(n-m,m)$. Every element $U\in U(n-m,m)$ can be represented as $U=e^{i\psi}\cdot V$ with $\psi\in\RR$ and $V\in SU(n-m,m)$. Note that there are no non-trivial homomorphisms from the unit circle into $\RR_{+}$ since $\RR_{+}$ has no non-trivial compact subgroups. Also, there are no non-trivial homomorphisms from $SU(n-m,m)$ into $\RR_{+}$ since the kernel of any such homomorphism is a proper normal subgroup of $SU(n-m,m)$ of positive dimension, and $SU(n-m,m)$ is a simple group. Thus, $\Lambda$ is constant on $U(n-m,m)$ and hence on all of $G_0(M)$. It then follows that, in coordinates in which $\hbox{Aut}_0(M)$ is linear, the function $F$ is invariant under all linear transformations of the $z$-variables from $U(n-m,m)$ and therefore depends only on $\langle z, z\rangle$ and $u$. Conditions (\ref{formnfcond}) imply that $F_{2\overline{2}}\equiv 0$, $F_{3\overline{3}}\equiv 0$. Thus, $F$ has the form
\begin{equation}
F(z,\overline{z},u)=\displaystyle\sum_{k=4}^{\infty}C_k(u)\langle z, z\rangle^k,\label{formnsq}
\end{equation}
where $C_k(u)$ are real-valued analytic functions of $u$, and for some $k$ we have $C_k(u)\not\equiv 0$. Note, in particular, that if $d_0(M)=n^2$, then 0 is an umbilic point of $M$. 

Conversely, if $M$ is given in the normal form by an equation
$$
v=\langle z,z\rangle+F(z,\overline{z},u),
$$
with $F\not\equiv 0$ of the form (\ref{formnsq}), then $\hbox{Aut}_0(M)$ contains all linear transformations (\ref{formlinearform}) with 
$U\in U(n-m,m)$, $\lambda=1$ and $\sigma=1$, and therefore $d_0(M)=n^2$. For $n>2m$ and for $n=2m$ with $G_0(M)=U(m,m)$, $\hbox{Aut}_0(M)$ clearly coincides with the group of all transformations of the form   
\begin{equation}
\begin{array}{lll}
z&\mapsto& Uz,\\
w&\mapsto& w.
\end{array}\label{linear2}
\end{equation}
where $U\in U(n-m,m)$. If $n=2m$ and $G_0(M)=U'(m,m)$, then $\hbox{Aut}_0(M)$ consists of all mappings  
$$
\begin{array}{lll}
z&\mapsto& Uz,\\
w&\mapsto &\sigma w,
\end{array}
$$
where $U\in U'(m,m)$, $\langle Uz,Uz\rangle=\sigma\langle z,z\rangle$, $\sigma=\pm 1$.

We will now concentrate on the case $0<d_0(M)<n^2$ (hence assuming that $n\ge 2$). For the strongly pseudoconvex case we obtain the following

\begin{theorem}\label{theorm0}\sl Let $M$ be a strongly pseudoconvex real-analytic non-spherical hypersurface in $\CC^{n+1}$ with $n\ge 2$ (here $m=0$) given in normal coordinates in which $\hbox{Aut}_0(M)$ is linear. Then the following holds
\smallskip\\

\noindent (i) $d_0(M)\ge n^2-2n+3$ implies $d_0(M)=n^2$;
\smallskip\\

\noindent (ii) if $d_0(M)=n^2-2n+2$, after a linear change of the $z$-coordinates the equation of $M$ takes the form
\begin{equation}
v=\sum_{\alpha=1}^n|z_{\alpha}|^2+F(z,\overline{z},u),\label{formm}
\end{equation}
where $F$ is a function of $|z_1|^2$, $\langle z,z\rangle:=\displaystyle\sum_{\alpha=1}^n|z_{\alpha}|^2$  and $u$:
\begin{equation}
F(z,\overline{z},u)=\sum_{p+q\ge 4}C_{pq}(u)|z_1|^{2p}\langle z,z\rangle^q,\label{formlessnsq}
\end{equation}
where $C_{pq}(u)$ are real-valued analytic functions of $u$, and $C_{pq}(u)\not\equiv 0$ for some $p,q$ with $p>0$;
\smallskip\\

\noindent (iii) if a hypersurface $M$ is given in the form described in (ii) (without assuming the linearity of $\hbox{Aut}_0(M)$ a priori), the group $\hbox{Aut}_0(M)$ coincides with the group of all mappings of the form (\ref{linear2}), where $U\in U(1)\times U(n-1)$ (with $U(1)\times U(n-1)$ realized as a group of block-diagonal matrices in the standard way).
\end{theorem}

\begin{corollary}\label{cor1}\sl If $M$ is a strongly pseudoconvex real-analytic hypersurface in $\CC^{n+1}$, and the dimension of $\hbox{Aut}_0(M)$ is greater than or equal to $n^2-2n+2$, then the origin is an umbilic point of $M$.
\end{corollary}

For the case $m\ge1$ we prove the following

\begin{theorem}\label{theorm1}\sl Let $M$ be a Levi non-degenerate real-analytic non-spherical hypersurface in $\CC^{n+1}$ with $m\ge 1$. Then the following holds
\smallskip\\

\noindent (i) $d_0(M)\ge n^2-2n+4$ implies $d_0(M)=n^2$;
\smallskip\\

\noindent (ii) if $d_0(M)=n^2-2n+3$, the group $\hbox{Aut}_0(M)$ is linearizable and in some normal coordinates in which $\hbox{Aut}_0(M)$ is linear, the equation of $M$ takes the form
\begin{equation}
\begin{array}{lll}
v&=&\displaystyle2\hbox{Re}\,z_1\overline{z}_n+2\hbox{Re}\,z_2\overline{z}_{n-1}+\dots+
2\hbox{Re}\,z_m\overline{z}_{n-m+1}+\\
&&\displaystyle\sum_{\alpha=m+1}^{n-m}|z_{\alpha}|^2
+F(z,\overline{z},u),
\end{array}\label{formm1minus}
\end{equation}
where $F$ is a function of $|z_n|^2$, $\langle z,z\rangle:=\displaystyle 
2\hbox{Re}\,z_1\overline{z}_n+2\hbox{Re}\,z_2\overline{z}_{n-1}+\dots+
2\hbox{Re}\,z_m\overline{z}_{n-m+1}+\sum_{\alpha=m+1}^{n-m}|z_{\alpha}|^2$ and $u$:
\begin{equation}
F(z,\overline{z},u)=\sum C_{rpq}u^r|z_n|^{2p}\langle z,z\rangle^q,\label{formlessnsq1}
\end{equation}
where at least one of $C_{rpq}\in\RR$ is non-zero, the summation is taken over $p\ge 1$, $q\ge 0$, $r\ge 0$ such that $(r+q-1)/p=s$ with $s\ge -1/2$ being a fixed rational number, and 
$$
F(z,\overline{z},u)=\displaystyle\sum_{k,\overline{l}\ge 2}F_{k\overline{l}}(z,\overline{z},u),
$$
where $F_{2\overline{3}}=0$ and identities (\ref{formnfcond}) hold for $F_{2\overline{2}}$ and $F_{3\overline{3}}$; 
\smallskip\\

\noindent (iii) if a hypersurface $M$ is given in the form described in (ii) (without assuming the linearity of $\hbox{Aut}_0(M)$ a priori), the group $\hbox{Aut}_0(M)$ coincides with the group of all mappings of the form
\begin{equation}
\begin{array}{lll}
z&\mapsto& |\mu|^{1/(s+1)} Uz,\\
w&\mapsto& |\mu|^{2/(s+1)} w,
\end{array}\label{mapwithlambda2}
\end{equation} 
with $U\in S$, where $S$ is the group introduced in Lemma \ref{alglemma} below, and $\mu$ is a parameter in this group (see formula (\ref{formalgebra})). 
\end{theorem}

\begin{corollary}\label{cor2}\sl Let $M$ be a Levi non-degenerate real-analytic hypersurface in $\CC^{n+1}$, with $n\ge 2$ and $m\ge 1$, and assume that the dimension of $\hbox{Aut}_0(M)$ is greater than or equal to $n^2-2n+3$. If the origin is a non-umbilic point of $M$, then in some normal coordinates the equation of $M$ takes the form
\begin{equation}
\begin{array}{lll}
v&=&2\hbox{Re}\,z_1\overline{z}_n+2\hbox{Re}\,z_2\overline{z}_{n-1}+\dots+\\
&&\displaystyle2\hbox{Re}\,z_m\overline{z}_{n-m+1}+\sum_{\alpha=m+1}^{n-m}|z_{\alpha}|^2\pm |z_n|^4.
\end{array}\label{formmx}
\end{equation}
\end{corollary}

We remark that hypersurfaces (\ref{formmx}) occur in \cite{P} in connection with studying unbounded homogeneous domains in complex space.

The proofs of Theorems \ref{theorm0} and \ref{theorm1} are given in Sections \ref{strong} and \ref{m1} respectively. Before proceeding we wish to acknowledge that this work was initiated while the first author was visiting the Mathematical Sciences Institute of the Australian National University.

\section{The Strongly Pseudoconvex Case}\label{strong}
\setcounter{equation}{0}

First of all, we note that in this case the mapping $\Lambda$ defined in (\ref{formlambda}) is constant, that is, $\lambda_{\varphi}=1$ for all $\varphi\in\hbox{Aut}_0(M)$. This follows from the fact that all eigenvalues of $U_{\varphi}$ are unimodular, or, alternatively, from the compactness of $G_0(M)$ and the observation that $\RR_{+}$ does not have non-trivial compact subgroups. Next, by a linear change of the $z$-coordinates the matrix $H$ can be transformed into the identity matrix $E_n$, and for the remainder of this section we assume that $H=E_n$. Hence we assume that the equation of $M$ is written in the form (\ref{formm}), where the function $F$ satisfies the normal form conditions.   

It is shown in Lemma 2.1 of \cite{IK} that any closed connected subgroup of the unitary group $U(n)$ of dimension $n^2-2n+3$ or larger is either $SU(n)$ or $U(n)$ itself. Hence, if $d_0(M)\ge n^2-2n+3$, we have $G_0(M)\supset SU(n)$, and therefore $F(z,\overline{z},u)$ is invariant under all linear transformations of the $z$-variables from $SU(n)$. This implies that $F(z,\overline{z},u)$ is a function of $\langle z,z\rangle$ and $u$, which gives that $F(z,\overline{z},u)$ is invariant under the action of the full unitary group $U(n)$ and thus $d_0(M)=n^2$, as stated in (i).

The proof of part (ii) of the theorem is also based on Lemma 2.1 of \cite{IK}. For the case $d_0(M)=n^2-2n+2$ the lemma gives that the connected identity component $G_0^c$ of $G_0$ is either conjugate in $U(n)$ to the subgroup $U(1)\times U(n-1)$ realized as block-diagonal 
matrices, or, for $n=4$, contains a subgroup conjugate in $U(n)$ to $Sp_{2,0}$. If the latter is the case, then, since $Sp_{2,0}$ acts transitively on the sphere of dimension 7 in $\CC^4$, $F(z,\overline{z},u)$ is a function of $\langle z,z\rangle$ and $u$, which implies that $F(z,\overline{z},u)$ is invariant under the action of the full unitary group $U(4)$ and thus $d_0(M)=16$, which is impossible. Hence $G_0^c$ is conjugate to $U(1)\times U(n-1)$, and therefore, after a unitary change of the $z$-coordinates, the equation of $M$ can be written in the form (\ref{formm}) where the function $F$ depends on $|z_1|^2$, $\langle z,z\rangle':=\displaystyle\sum_{\alpha=2}^n|z_{\alpha}|^2$ and $u$. Clearly, $\langle z,z\rangle'=\langle z,z\rangle-|z_1|^2$, and $F$ can be written as a function of $|z_1|^2$, $\langle z,z\rangle$ and $u$ as in (\ref{formlessnsq}). Next, conditions (\ref{formnfcond}) imply that $F_{2\overline{2}}\equiv 0$, $F_{3\overline{3}}\equiv 0$, and thus the summation in (\ref{formlessnsq}) is taken over $p,q$ such that $p+q\ge 4$. Further, if $C_{pq}\equiv 0$ for all $p>0$, $F$ has the form (\ref{formnsq}) and therefore $G_0=U(n)$ which is impossible. Thus for some $p,q$ with $p>0$ we have $C_{pq}\not\equiv 0$, and (ii) is established.         

If $M$ is given in the normal form and is written as in (\ref{formm}), (\ref{formlessnsq}), $\hbox{Aut}_0(M)$ clearly contains all maps of the form (\ref{linear2}) with  $U\in U(1)\times U(n-1)$. Hence $d_0(M)\ge n^2-2n+2$. If $d_0(M)> n^2-2n+2$, then by part (i) of the theorem, $d_0(M)=n^2$ and hence $G_0(M)=U(n)$. Then $F$ has the form (\ref{formnsq}) which is impossible because for some $p,q$ with $p>0$ the function $C_{pq}$ does not vanish identically. Thus $d_0(M)=n^2-2n+2$, and hence $G_0^c(M)=U(1)\times U(n-1)$. It is not hard to show that $G_0(M)$ is connected (note, for example, that by an argument given in the introduction, $\hbox{Aut}_0(M)$ is linear in these coordinates), and therefore $\hbox{Aut}_0(M)$ coincides with the group of all mappings of the form (\ref{linear2}), where $U\in U(1)\times U(n-1)$.    

Thus, (iii) is established, and the theorem is proved.\qed

\section{The Case of $m\ge1$}\label{m1}
\setcounter{equation}{0}

We start with the following algebraic lemma.

\begin{lemma}\label{alglemma}\sl Let $G\subset U(n-m,m)$ be a connected real algebraic subgroup of $GL_n(\CC)$, $n\ge 2m$, $m\ge 1$, with Hermitian form preserved by $U(n-m,m)$ written as
\begin{equation}
\left(
\begin{array}{lllllll}
&&&&&&1\\
&\hbox{\Huge 0}&&&&\ddots&\\
&&&&1&&\\
&&&E_{n-2m}&&&\\
&&1&&&&\\
&\ddots&&&&\hbox{\Huge 0}&\\
1&&&&&&
\end{array}
\right),\label{formformaminus}
\end{equation}
where $E_{n-2m}$ is the $(n-2m)\times(n-2m)$ identity matrix, and the number of 1's on each side of $E_{n-2m}$ is $m$. Then the following holds
\smallskip\\

\noindent (a) if $\hbox{dim}\,G\ge n^2-2n+4$, we have either $G=SU(n-m,m)$, or $G=U(n-m,m)$;
\smallskip\\

\noindent (b) if $\hbox{dim}\,G=n^2-2n+3$, the group $G$ either is conjugate in $U(n-m,m)$ to the group $S$ that consists of all matrices of the form
\begin{equation}
\left(
\begin{array}{ccc}
\mu & -\mu\overline{x}^TH'A & c\\
\vspace{0.05cm} &&\\
0 & A & x\\
\vspace{0.05cm} &&\\
0 & 0 & 1/\overline{\mu}
\end{array}
\right),\label{formalgebra}
\end{equation}
where $\mu,c\in\CC$, $\mu\ne 0$, $x\in\CC^{n-2}$, $A\in U(n-m-1,m-1)$ (i.e., $A$ is an $(n-2)\times(n-2)$-matrix with complex elements such that $A^TH'\overline{A}=H'$ with $H'$ obtained from matrix (\ref{formformaminus}) by removing the first and the last columns and rows), and the following holds
$$
2\hbox{Re}\,\displaystyle\frac{c}{\mu}+x^TH'\overline{x}=0,
$$
or, if $n=4$ and $m=2$, coincides with $e^{i\RR}\Bigl(Sp_4(B,\CC)\cap SU(2,2)\Bigr)$, or, if $n=2$ and $m=1$, coincides with $SU(1,1)$. Here the subgroup $Sp_4(B,\CC)\subset GL_4(\CC)$ consists of matrices preserving a non-degenerate skew-symmetric bilinear form $B$ equivalent to the form given by the matrix
\begin{equation}
B_0:=\left(
\begin{array}{ll}
0 & E_2\\
-E_2 & 0
\end{array}
\right),\label{b0}
\end{equation}
where $E_2$ is the $2\times 2$ identity matrix,
\end{lemma}

\noindent {\bf Proof:} Let $V\subset U(n-m,m)$ be a real algebraic subgroup of $GL_n(\CC)$ such that $\hbox{dim}\, V\ge n^2-2n+3$. Consider $V_1:=V\cap SU(n-m,m)$. Clearly, $\hbox{dim}\, V_1\ge n^2-2n+2$. Let $V^{\CC}_1\subset SL_n(\CC)$ be the complexification of $V_1$. We have $\hbox{dim}_{\CC}V^{\CC}_1\ge n^2-2n+2$. Consider the maximal complex closed subgroup $W(V)\subset  SL_n(\CC)$ that contains $V^{\CC}_1$. Clearly, $\hbox{dim}_{\CC} W(V)\ge n^2-2n+2$. All closed maximal subgroups of $SL_n(\CC)$ had been classified (see \cite{D}), and the lower bound on the dimension of $W(V)$ gives that either $W(V)=SL_n(\CC)$, or $W(V)$ is conjugate to one of the parabolic subgroups
$$
\begin{array}{l}
P^1:=\left\{
\left(
\begin{array}{cc}
1/\det C & b\\
0 & C
\end{array}
\right),\, b\in\CC^{n-1},\, C\in GL_{n-1}(\CC)\right\},\\
\vspace{0.2cm}\\
P^2:=\left\{
\left(
\begin{array}{cc}
C & b\\
0 & 1/\det C
\end{array}
\right),\,b\in\CC^{n-1},\, C\in GL_{n-1}(\CC)\right\}
\end{array}
$$
(note that $P^1=P^2$ for $n=2$), or, for $n=4$, $W(V)$ is conjugate to $Sp_4(\CC)$.

Suppose that for some $g\in SL_n(\CC)$ and $j\in\{1,2\}$ we have $g^{-1}W(V)g=P^j$. It is not hard to show that, due to the lower bound on the dimension of $W(V)$, $g$ can be chosen to belong to $SU(n-m,m)$. Then $g^{-1}V_1g\subset P^j\cap SU(n-m,m)$. It is easy to compute the intersections $P^j\cap SU(n-m,m)$ for $j=1,2$ and see that they are equal and coincide with the group $S_1$ of matrices of the form (\ref{formalgebra}) with determinant 1. Clearly, $\hbox{dim}\,S_1=n^2-2n+2\le\hbox{dim}\,V_1$ and therefore $V_1$ is conjugate to $S_1$ in $SU(n-m,m)$.

Suppose now that $n=4$ and for some $g\in SL_4(\CC)$ we have $g^{-1}W(V)g=Sp_4(\CC)$. In particular, $g^{-1}V_1g\subset Sp_4(\CC)\cap g^{-1}SU(4-m,m)g$ (here we have either $m=1$, or $m=2$). It can be shown that $\hbox{dim}\, Sp_4(\CC)\cap g^{-1}SU(3,1)g\le 6$ for all $g\in SL_4(\CC)$. At the same time we have $\hbox{dim}\, V_1\ge 10$. Hence $W(V)$ in fact cannot be conjugate to $Sp_4(\CC)$, if $m=1$. Therefore, $m=2$, and $V_1\subset g Sp_4(\CC)g^{-1}\cap SU(2,2)=Sp_4(B,\CC)\cap SU(2,2)$, where $B$ is some non-degenerate skew-symmetric bilinear form. It is straightforward to show that $Sp_4(B,\CC)\cap SU(2,2)$ is connected and $\hbox{dim}\, Sp_4(B,\CC)\cap SU(2,2)\le 10$. Therefore $V_1=Sp_4(B,\CC)\cap SU(2,2)$. 

Suppose now that $\hbox{dim}\,G\ge n^2-2n+4$. Then $\hbox{dim}\,G_1\ge n^2-2n+3$, and the above considerations give that $W(G)=SL_n(\CC)$. Hence $G_1=SU(n-m,m)$ which implies that either $G=SU(n-m,m)$, or $G=U(n-m,m)$, thus proving (a).

Let $\hbox{dim}\,G=n^2-2n+3$. In this case we have either $\hbox{dim}\,G_1=n^2-2n+2$, or $G=SU(1,1)$, if $n=2$, $m=1$. In the first case we obtain that $G_1$ either is conjugate to $S_1$ in $SU(n-m,m)$, or, for $n=4$ and $m=2$ coincides with $Sp_4(B,\CC)\cap SU(2,2)$ for some non-degenerate skew-symmetric bilinear form $B$ equivalent to the form $B_0$ defined in (\ref{b0}). This gives that $G$ in the first case either is conjugate to $S$ in $U(n-m,m)$, or for $n=4$ and $m=2$ coincides with $e^{i\RR}\Bigl(Sp_4(B,\CC)\cap SU(2,2)\Bigr)$, and (b) is established. 

The lemma is proved.\qed
\medskip\\

We will now prove Theorem \ref{theorm1}. Suppose first that $d_0(M)\ge n^2-2n+4$ and assume that $H$ is written in the diagonal form with 1's in the first $n-m$ positions and  $-1$'s in the last $m$ positions on the diagonal. Lemma \ref{alglemma} gives that either $G_0^c(M)=SU(n-m,m)$ (in which case $n\ge 3$), or $G_0(M)=U(n-m,m)$, or, for $n=2m$, $G_0(M)=U'(m,m)$. If $G_0(M)\supset U(n-m,m)$, then $d_0(M)=n^2$, and (i) is established. Assume that $G_0(M)\supset SU(n-m,m)$. Suppose that $m\ge 2$. Then $G_0$ contains the product $R:=SU(n-m)\times SU(m)$ realized as block-diagonal matrices. Arguing as in the introduction, we obtain that in some normal coordinates all elements of the compact group $\hat R:=\Phi^{-1}(R)$ can be written in the form (\ref{linear2}) and thus $F$ is a function of $\langle z,z\rangle_{+}:=\sum_{j=1}^{n-m}|z_j|^2$, $\langle z,z\rangle_{-}:=\sum_{j=n-m+1}^{n}|z_j|^2$, and $u$. Hence all elements of odd weight in the weight decomposition for $F$ are zero. This shows that $F_{\gamma+1}\equiv 0$, and identity (\ref{eq3}) again implies that $\hbox{Aut}_0(M)$ becomes linear after a change of coordinates of the form (\ref{fractionallintransfm}). If $m=1$, $\hbox{Aut}_0(M)$ is linearizable by \cite{Ezh1}, \cite{Ezh2}.   

Therefore, there exist normal coordinates where the corresponding function $F$ is invariant under all linear transformations of the $z$-variables from $SU(n-m,m)$. This implies that $F$ is in fact invariant under all linear transformations of the $z$-variables from $U(n-m,m)$. Hence $d_0(M)=n^2$, and (i) is established.       

Suppose now that $d_0(M)=n^2-2n+3$. By a linear change of the $z$-coordinates the matrix $H$ can be transformed into matrix (\ref{formformaminus}), and from now on we assume that $H$ is given in this form. Hence the equation of $M$ is written as in (\ref{formm1minus}), where the function $F$ satisfies the normal form conditions. Arguing as in the preceding paragraph, we see that for $n=2$, $m=1$, the group $G_0^c$ cannot coincide with $SU(1,1)$. Assume first that after a linear change of the $z$-coordinates preserving the form $H$ the group $G_0^c(M)$ coincides with $S$. Then $G_0(M)$ contains the compact subgroup $Q=\{e^{it}\cdot E_n,\,t\in\RR\}$, where $E_n$ is the $n\times n$ identity matrix.
The argument based on identity (\ref{eq3}) that we gave in the introduction, again yields that $\hbox{Aut}_0(M)$ is linearizable. Passing to coordinates in which $\hbox{Aut}_0(M)$ is linear, we obtain that for every $U\in S$ the equation of $M$ is invariant under the linear transformation
\begin{equation}
\begin{array}{lll}
z&\mapsto& \lambda_U Uz,\\
w&\mapsto& \lambda_U^2 w,
\end{array}\label{mapwithlambda}
\end{equation}
where $\lambda_U=\Lambda(U)$. The group $S$ contains $U(n-m-1,m-1)$ realized as the subgroup of all matrices of the form (\ref{formalgebra}) with $\mu=1$, $c=0$, $x=0$. Since $\Lambda$ is constant on $U(n-m-1,m-1)$, we have $\lambda_U=1$ for all $U\in U(n-m-1,m-1)$. Therefore, the function $F(z,\overline{z},u)$ depends on $z_1$, $z_n$, $\overline{z}_1$, $\overline{z}_n$, $\langle z,z\rangle':=\displaystyle 2\hbox{Re}\,z_2\overline{z}_{n-1}+\dots+
2\hbox{Re}\,z_m\overline{z}_{n-m+1}+\sum_{\alpha=m+1}^{n-m}|z_{\alpha}|^2$ and $u$. Clearly, $\langle z,z\rangle'=\langle z,z\rangle-2\hbox{Re}\,z_1\overline{z_n}$, and $F$ can be written as follows
$$
F(z,\overline{z},u)=\sum_{r,q\ge 0} D_{rq}(z_1,z_n,\overline{z}_1,\overline{z}_n)u^r\langle z,z\rangle^q,
$$
where $D_{rq}$ are real-analytic.

We will now determine the form of the functions $D_{rq}$. The group $S$ contains the subgroup $I$ of all matrices as in (\ref{formalgebra}) with $|\mu|=1$, $x=0$ and $A=E_{n-2}$, where $E_{n-2}$ is the $(n-2)\times (n-2)$ identity matrix. Since every eigenvalue of any $U\in I$ is unimodular, we have $\lambda_U=1$ for all $U\in I$, and therefore $D_{rq}$ is invariant under all linear transformations from $I$. It is straightforward to show (see also \cite{Ezh2}) that any polynomial of $z_1$, $z_n$, $\overline{z_1}$, $\overline{z_n}$ invariant under all linear transformations from $I$ is a function of $\hbox{Re}\,z_1\overline{z}_n$ and $|z_n|^2$, and hence every $D_{rq}$ has this property. Let further $J$ be the subgroup of $S$ given by the conditions $\mu=1$, $A=E_{n-2}$. For every $U\in J$ we also have $\lambda_U=1$, and hence $D_{rq}$ is invariant under all linear transformations from $J$. It is then easy to see that $D_{rq}$ has to be a function of $|z_n|^2$ alone. Thus, the function $F$ has the form (\ref{formlessnsq1}), and it remains to show that the summation in (\ref{formlessnsq1}) is taken over $p\ge 1$, $q\ge 0$, $r\ge 0$ such that $(r+q-1)/p=s$, where $s\ge -1/2$ is a fixed rational number.

Let $K$ be the 1-dimensional subgroup of $S$ given by the conditions $\mu>0$, $c=0$, $x=0$, $A=E_{n-2}$. It is straightforward to show that every homomorphism $\Psi:K\ra\RR_{+}$ has the form $U\mapsto \mu^{\alpha}$, where $\alpha\in\RR$. Considering $\Psi=\Lambda|_K$ we obtain that there exists $\alpha\in\RR$ such that for every $U\in K$ we have $\lambda_U=\mu^{\alpha}$. We will now prove that $\alpha\ne 0$. Indeed, otherwise $F$ would be invariant under all linear transformations from $K$ and therefore would be a function of $\langle z,z\rangle$ and $u$, which implies that $G_0(M)\supset U(n-m,m)$. This contradiction shows that $\alpha\ne 0$ and hence $\lambda_U\ne 1$ for every $U\in K$ with $\mu\ne 1$.

Plugging a mapping of the form (\ref{mapwithlambda}) with $U\in K$, $\mu\ne 1$, into equation (\ref{formm1minus}), where $F\not\equiv 0$ has the form (\ref{formlessnsq1}) we obtain that, if $C_{rpq}\ne 0$, then
\begin{equation}
\lambda_U^{r+p+q-1}=\mu^p.\label{rellammu}
\end{equation}
The equation of $M$ is written in the normal form, hence $p+q\ge 2$ and $r+p+q-1\ge 1$. Since $\lambda_U\ne 1$, we obtain that $p\ge 1$. Further, (\ref{rellammu}) implies
$$
\lambda_U^{(r+p+q-1)/p}=\mu,  
$$
and, since the right-hand side in the above identity does not depend on $r,p,q$, for all non-zero coefficients $C_{rpq}$ the ratio $(r+q-1)/p$ must have the same value; we denote it by $s$. Clearly, $s$ is a rational number and $s\ge -1/2$. We also remark that $\alpha=p/(r+p+q-1)=1/(s+1)$.

Assume now that $n=4$, $m=2$ and $G_0^c(M)$ coincides with $e^{i\RR}\Bigl(Sp_4(B,\CC)\cap SU(2,2)\Bigr)$ for some non-degenerate skew-symmetric non-degenerate bilinear form $B$ equivalent to the form $B_0$ defined in (\ref{b0}). Then $G_0(M)$ contains the compact subgroup $Q=\{e^{it}\cdot E_4,\,t\in\RR\}$, where $E_4$ is the $4\times 4$ identity matrix. Arguing as above, we obtain that $\hbox{Aut}_0(M)$ is linearizable. Further, it is straightforward to prove that $Sp_4(B,\CC)\cap SU(2,2)$ is a real form of $Sp_4(B,\CC)$ and therefore is simple. Hence there does not exist a non-trivial homomorphism from $Sp_4(B,\CC)\cap SU(2,2)$ into $\RR_{+}$. Further, since $\RR_{+}$ does not have non-trivial compact subgroups, any homomorphism from the unit circle into $\RR_{+}$ is constant. Hence $\Lambda$ is constant on $G_0(M)$. This implies that $F$ is invariant under all linear transformations from $Sp_4(B,\CC)\cap SU(2,2)$. It can be shown that this group acts transitively on any pseudosphere in $\CC^4$ given by the equation $\langle z,z\rangle =r$, which yields that $F$ is a function of  $\langle z,z\rangle$ and $u$ and hence $d_0(M)=n^2$. This contradiction proves that in fact $G_0^c(M)\ne e^{i\RR}\Bigl(Sp_4(B,\CC)\cap SU(2,2)\Bigr)$ for $n=4$, $m=2$. Thus, (ii) is established.

Suppose that $M$ is given in the normal form, written as in (\ref{formm1minus}), (\ref{formlessnsq1}), and the summation in (\ref{formlessnsq1}) is taken over $p\ge 1$, $q\ge 0$, $r\ge 0$ such that $(r+q-1)/p=s$, where $s\ge -1/2$ is a fixed rational number. 
Set $\alpha=1/(s+1)$ and for every $U\in S$ define $\lambda_U=|\mu|^{\alpha}$. It is then straightforward to verify that every mapping of the form (\ref{mapwithlambda}) with $U\in S$ is an automorphism of $M$. Therefore, $G_0(M)$ contains $S$ and hence $d_0(M)\ge n^2-2n+3$. If $d_0(M)>n^2-2n+3$, then by part (i) of the theorem, $d_0(M)=n^2$ and hence $G_0(M)\supset U(n-m,m)$. Then $F$ is a function of $\langle z,z\rangle$ and $u$, which is impossible since for every non-zero $C_{rpq}$ we have $p\ge 1$. Hence $d_0(M)= n^2-2n+3$ and hence $G_0^c(M)=S$. Finally, observe that by an argument given in the introduction, $\hbox{Aut}_0(M)$ is linear in these coordinates. It is now straightforward to show that $\hbox{Aut}_0(M)$ coincides with the group of all mappings of the form (\ref{mapwithlambda2}).

Thus, (iii) is established, and the theorem is proved.\qed

{\obeylines
School of Mathematics and Statistics
University of South Australia
Mawson Lakes Blvd
Mawson Lakes
South Australia 5091
AUSTRALIA
E-mail: vladimir.ejov@unisa.edu.au
\hbox{ \ \ }

\hbox{ \ \ }
Department of Mathematics
The Australian National University
Canberra, ACT 0200
AUSTRALIA
E-mail: alexander.isaev@maths.anu.edu.au}

\end{document}